\documentclass[11pt,a4paper]{article}

\usepackage{amsfonts, amsmath, wasysym}
\usepackage{stmaryrd} 

\usepackage{tikz}
\usetikzlibrary{shapes, positioning}

\usepackage{amssymb,amsthm,
paralist
}

\usepackage{
latexsym,
}

\usepackage{url}

\definecolor{darkgreen}{rgb}{0,0.5,0}
\definecolor{darkred}{rgb}{0.7,0,0}
\usepackage[colorlinks, 
citecolor=darkgreen, linkcolor=darkred
]{hyperref}

\theoremstyle{plain}

\numberwithin{equation}{section}

\newcommand{\ca}{\ensuremath{{\mathcal A}}}

\newcommand{\cf}{\ensuremath{{\mathcal F}}}

\newcommand{\zb}{\bar{z}}

\newcommand{\pl}[2]{{\frac{\partial #1}{\partial #2}}}

\newcommand{\ga}{\gamma}

\newcommand{\ka}{\kappa}

\newcommand{\si}{\sigma}

\newcommand{\Si}{\Sigma}

\newcommand{\vph}{\varphi}
\newcommand{\ep}{\varepsilon}

\newcommand{\R}{\ensuremath{{\mathbb R}}}

\newcommand{\Z}{\ensuremath{{\mathbb Z}}}
\newcommand{\C}{\ensuremath{{\mathbb C}}}

\newcommand{\lap}{\Delta}

\newcommand{\grad}{\nabla}

\newcommand{\beq}{\begin{equation}}
\newcommand{\eeq}{\end{equation}}
\newcommand{\beqa}{\begin{equation}\begin{aligned}}
\newcommand{\eeqa}{\end{aligned}\end{equation}}
\newcommand{\brmk}{\begin{rmk}}
\newcommand{\ermk}{\end{rmk}}
\newcommand{\partref}[1]{\hbox{(\csname @roman\endcsname{\ref{#1}})}}
\newcommand{\half}{\frac{1}{2}}

\newcommand{\Hess}{{\mathrm{Hess}}}

\usepackage{soul}

 \newtheorem{thm}{Theorem}[section]

\newtheorem{lem}[thm]{Lemma}

\newtheorem{rmk}[thm]{Remark}

\title{\sc monotonicity of the modulus under \\ curve shortening flow}
\author{Arjun Sobnack and Peter M. Topping}
\date{14 April 2026}

\begin{document}

%

\parskip 8pt
\parindent 0pt

\maketitle

\begin{abstract}
Given two disjoint nested embedded closed curves in the plane, both evolving under curve shortening flow, we show that the modulus of the enclosed annulus is monotonically increasing in time. 
An analogous result holds within any ambient surface 
satisfying a lower curvature bound.
\end{abstract}

\section{Introduction}

The curve shortening flow evolves embedded curves in the direction of their geodesic curvature vector. In its simplest form we could start with a smooth embedded loop
in the plane, i.e.~an embedding $\ga:S^1\to\R^2$, and evolve it in time under the PDE
$$\pl{\ga}{t}=\vec\ka,$$
where the geodesic curvature vector $\vec\ka$ of the curve at a given time could be written with respect to a local arc-length parameter $s$ as $\vec\ka(s)=\ga''(s)$.
If we let $\nu$ be the outward unit normal and define the (scalar) geodesic curvature 
$\ka$ by $\vec\ka=-\ka\nu$, then 
the area $A$ enclosed by the evolving loop is governed by 
$$\frac{dA}{dt}=-\int \ka\, ds=-2\pi,$$
which limits its possible existence time to $\frac{A_0}{2\pi}$, where $A_0$ is the initial enclosed area.
Following earlier work of Gage and Hamilton \cite{gage_hamilton}, Grayson \cite{grayson} proved that 
a solution to this flow does indeed exist until time $T=\frac{A_0}{2\pi}$, with the curve looking like a tiny shrinking circle shortly before the final time.

This flow generalises in a straightforward way to more general ambient Riemannian manifolds than $\R^2$. In particular we are interested in evolving curves within 
a more general  Riemannian surface (see e.g.~\cite{grayson2}).
It is also natural to consider the evolution of curves other than loops. 
In this case we generally consider proper embeddings $\ga:\R\to\R^2$ and evolve them under the same equation. The most famous instance of such noncompact solutions is given by the so-called \textit{grim reaper}. This is a proper embedded curve $\ga_0:\R\to\R^2$ 
that when translated at unit speed in a given direction yields a solution to curve shortening flow, modulo reparametrisation. The standard approach to finding this curve is to write it as a graph of a function $x=\vph(y)$ and derive an ODE for $\vph$ that expresses that the curve translates in the $x$ direction under curve shortening flow. Standard ODE methods then yield that
$\vph(y)=-\log \cos y$ for $y\in (-\frac{\pi}2,\frac{\pi}2)$.

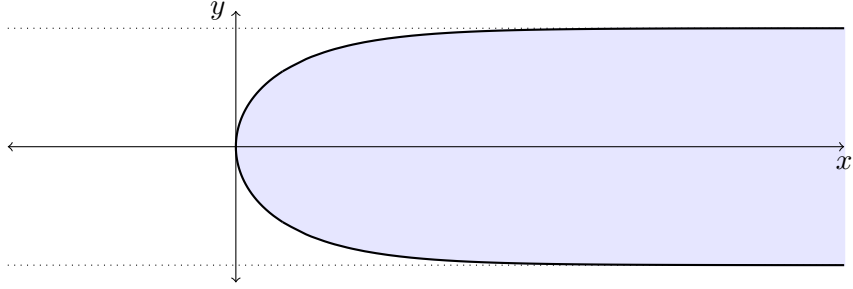
\begin{figure}
\centering
\begin{tikzpicture}

\draw[fill, scale=1, thick, domain=-1.1:1.1, smooth, variable=\y, blue!10]  plot ({-ln(cos(deg(\y))}, {\y}) -- (8,0);

\draw[fill, scale=1, thick, domain=-ln(cos(deg(1))):8, smooth, variable=\x, blue!10]  
plot ({\x},{rad(acos(exp(-\x))}) -- (8,0);

\draw[fill, scale=1, thick, domain=-ln(cos(deg(1))):8, smooth, variable=\x, blue!10]  
plot ({\x},{-rad(acos(exp(-\x))})  -- (8,0);

\draw[scale=1, thick, domain=-1:1, smooth, variable=\y, black]  plot ({-ln(cos(deg(\y))}, {\y});

\draw[scale=1, thick, domain=-ln(cos(deg(1))):8, smooth, variable=\x, black]  
plot ({\x},{rad(acos(exp(-\x))});

\draw[scale=1, thick, domain=-ln(cos(deg(1))):8, smooth, variable=\x, black]  
plot ({\x},{-rad(acos(exp(-\x))});

\draw [<->] (-3,0) -- (8,0) node[below]{$x$} ;
\draw [<->] (0,-1.8) -- (0,1.8) node[left]{$y$} ;

\draw [dotted] (-3,pi/2) -- (8,pi/2);
\draw [dotted] (-3,-pi/2) -- (8,-pi/2);

\end{tikzpicture}
\caption{Grim reaper}
\label{fig:grim}
\end{figure}

In this paper we take a different viewpoint on curve shortening flow and see it in terms of the geometry of the region it bounds. For example, the distinguishing feature of the 
grim reaper curve above turns out to be that the `smaller' region that it bounds,
shaded in Figure \ref{fig:grim}, is mapped to a half-space $\{z\in\C\ |\ \Re(z)> 1\}$
under the biholomorphic function 
$$z\mapsto e^z.$$
We are mainly concerned with the case that we have two disjoint nested embedded loops in the plane, or more generally two disjoint embedded loops in a Riemannian surface that bound an annulus, with both loops evolving under curve shortening flow. The avoidance principle, i.e.~the maximum principle, tells us that these two curves remain disjoint under the flow and therefore continue to bound an annulus. 
By virtue of being a Riemann surface with the topology of an annulus, this region  is biholomorphic to a cylinder.
While the boundary curves are nondegenerate this 
cylinder can be taken to be $(0,h)\times S^1$ for some $h\in (0,\infty)$, where $h$ is independent of the biholomorphism, and where we adopt the convention that 
$S^1:=\R/\Z$. Because we are considering circles $S^1$ of length $1$, the \textbf{modulus} of the annulus 
can be defined to be $h$.

In this note we demonstrate that this modulus is monotonically increasing as the boundary curves flow. 

\begin{thm}
\label{main_thm_mini}
Given two disjoint nested embedded loops in $\R^2$ that evolve under curve shortening flow over some time interval $t\in (0,T)$, the modulus of the enclosed annulus is strictly increasing in $t$.
\end{thm}

This theorem is a special case of the following  result within 
more general ambient surfaces.

\begin{thm}
\label{main_thm}
Suppose $M$ is a Riemannian surface with Gauss curvature $K$ bounded below
by $K_0\in \R$.
Given two disjoint embedded loops in $M$ that 
evolve under curve shortening flow over some time interval\, $t\in (0,T)$, and which
bound an annulus that evolves continuously in the Hausdorff sense,  the modulus 
$h$ 
of the enclosed annulus is controlled by
\beq
\label{logh_ineq}
\frac{d}{dt}\log h\geq K_0.
\eeq
If equality is attained in \eqref{logh_ineq} at some time
$t\in (0,T)$ then both $K_0=0$ and the enclosed annulus is isometric to a scaling of the flat cylinder 
$(0,h)\times S^1$.
\end{thm}

A particularly interesting ambient space here is the cylinder $\R\times S^1$ itself, with each curve a deformation of the circle $\{0\} \times S^1$.

The hypothesis that the annulus evolves continuously in the Hausdorff sense serves merely to prevent the annulus jumping in the case that two distinct annuli have the same boundary curves.

Although the boundary curves are evolving independently, the monotonicity of Theorem \ref{main_thm} is certainly not true in general if we 
evolve just one of the boundary curves. The two curves interact by 
being barriers for each other that evolve under the same equation. 
In general our motivation in this work is to 
learn how to prove regularity estimates and geometric constraints on a solution of curve shortening flow by virtue of it having weak proximity to a solution with known regularity. 
More concretely, in \cite{ST1} we proposed a framework
of which the following is a variant:

\begin{figure} 
\centering
\includegraphics[width=0.8\textwidth]{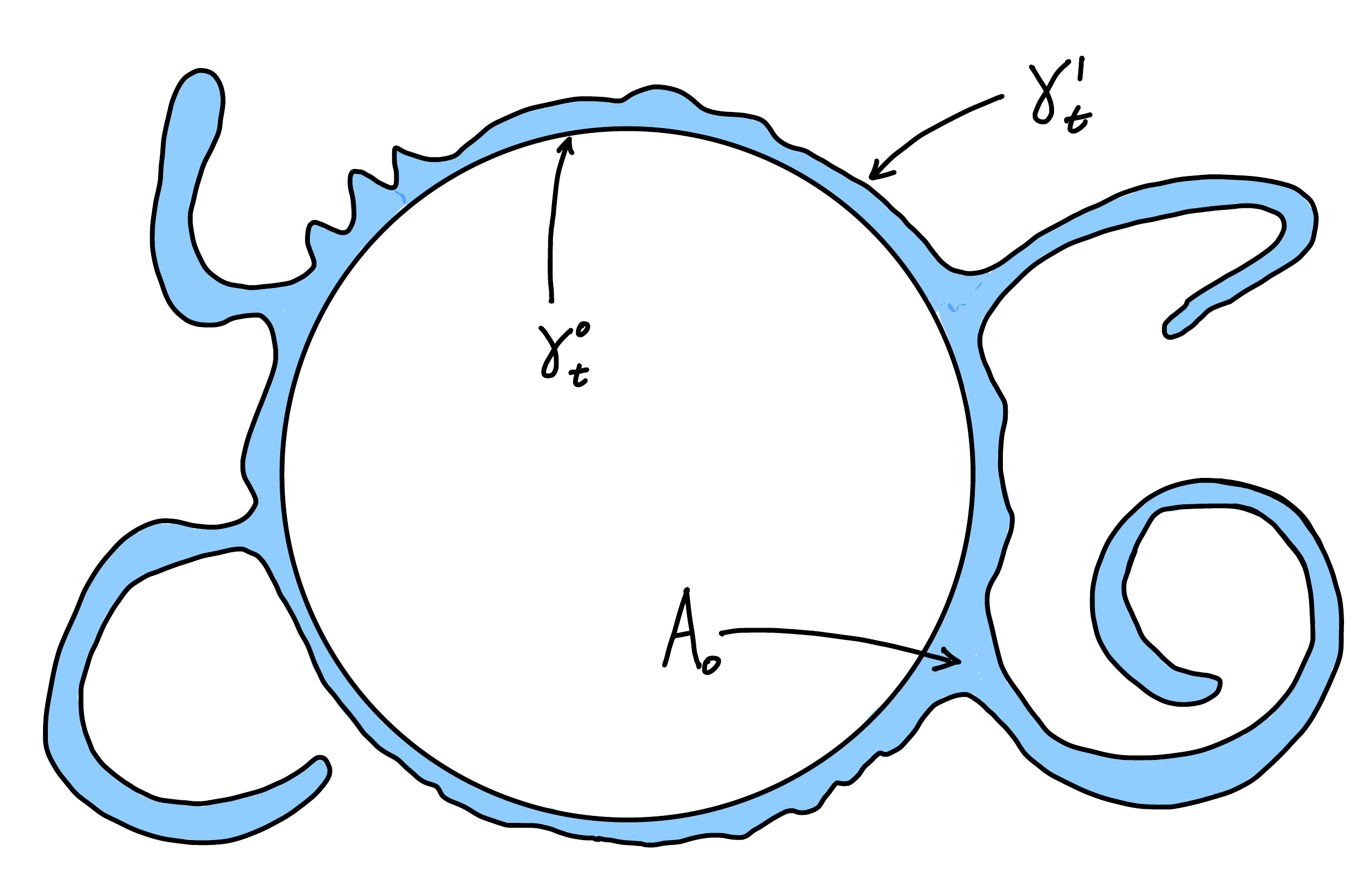}
\caption{Wild outer loop $\ga^1_t$ and regular inner loop $\ga^0_t$, enclosing a shaded annulus of area $A_0$.}
\label{A_0_fig}
\end{figure}

\begin{quote}
\it
Given two disjoint  embedded loops $\ga^0_t$ and $\ga^1_t$ in the plane, evolving under curve shortening flow for $t\in [0,T)$ and bounding an annulus of  area $A_0\ll T$, 
define the threshold time $\tau$ to be $\frac{A_0}{\pi}$.
Then after flowing beyond  the threshold time we expect to start being able to control the regularity of 
$\ga^1_t$ at time $\tau<t \ll T$ 
in terms of $t$, $A_0$ and the regularity of 
$\ga^0_t$ at time $t$. 
\end{quote}
Here when we talk of controlling the regularity we specifically would like to consider the length of the curve and bounds on its curvature and the derivatives of curvature.
An example would be when $\ga^0_t$ is a shrinking circle and 
$\ga^1_t$ is a loop outside this circle that starts with extremely wild regularity but nevertheless traps a small amount of area $A_0$ between $\ga^0_t$ and $\ga^1_t$ (with the trapped area necessarily constant in time). See Figure \ref{A_0_fig}.
We showed in \cite{ST1} that before that special time $\tau$ 
no uniform bounds on the regularity of $\ga^1_t$ in terms of $t$ and $A_0$ could hold.
Also in \cite{ST1}, we showed results exactly as proposed in the framework, in the case of graphical curve shortening flow. We believe that Theorem \ref{main_thm_mini} is likely to be a key ingredient in establishing the framework above in the non-graphical case such as in Figure \ref{A_0_fig}.

A bridge between this theory and the theory of Ricci flow is given by the 
work of the second author and H. Yin in their $L^1-L^\infty$ smoothing estimate from \cite{TY1}.

\emph{Acknowledgements:} 
AS was supported by EPSRC studentship EP/R513374/1.
PT  was supported by EPSRC grant EP/T019824/1.
For the purpose of open access, the authors have applied a Creative Commons Attribution (CC BY) licence to any author accepted manuscript version arising.

\section{The modulus of an annulus}

It is a consequence of the Uniformisation theorem that every Riemann surface $\Si$ that is homeomorphic to an annulus must either be biholomorphic to a
cylinder $(0,h)\times S^1$, or to one of 
$(0,\infty)\times S^1$ or $\R\times S^1$.

Because of our convention that $S^1:=\R/\Z$, i.e. $S^1$ has length $1$,
in the former cases we can define the \textit{modulus} of the Riemann surface to be 
$h$.
We then see that $h$ could be alternatively written in terms of the 
(conformally invariant) Dirichlet energy 
$$E(v)=i\int_\Si |v_z|^2 dz\wedge d\zb$$
of functions $v$ within the class $\cf$ of functions $v\in C^1(\Si)$ that when viewed as functions on the open cylinder $(0,h)\times S^1$
would extend continuously to the closed cylinder 
$[0,h]\times S^1$, taking values $0$ at one end and $1$ at the other end.
Indeed, if we define the capacity of $\Si$ to be
\beq
\label{var_prob}
\mathrm{Cap}(\Si):=\inf \{ E(v) : v\in \cf\},
\eeq
then this infimum is attained by $u:=x/h$, where $x$ is the coordinate along the cylinder,
which has energy
$$E(u)=i\int_{(0,h)\times S^1}{\textstyle |\frac1{2h}|^2}(-2i)dx dy=\frac{1}{2h}$$
since we are taking $S^1$ of length $1$.
We thus see that 
\beq
\label{var_prob2}
h=\frac{1}{2\,\mathrm{Cap}(\Si)}.
\eeq
Alternatively, we can equip $\Si$ with its unique complete conformal hyperbolic metric, which would be written 
$$\rho^2(dx^2+dy^2)\qquad\text{where}\qquad \rho(x)=\frac{\pi}{h\sin (\frac{\pi x}{h})}$$
on the cylinder, and then define the modulus in terms of the length $\ell$ of the shortest closed geodesic separating the two boundary components 
using the formula
$$h=\frac{\pi}{\ell}.$$
These geometric considerations make it straightforward to read off  properties of the modulus. An initial observation is that if we isotop $\Si$ to a smaller annulus strictly within $\Si$, then the modulus will strictly decrease.
We are primarily concerned in this section with deriving a formula for the rate of change of the modulus,
or equivalently capacity, of an annulus $\ca_t$ that is a  subset of a fixed ambient Riemannian surface $M$,
with $\ca_t$ evolving continuously in the Hausdorff sense, such that the boundary $\partial \ca_t$ consists of two disjoint smooth embedded loops in $M$ that vary 
smoothly in time for $t\in (-\ep,\ep)$.
If the curves vary in the direction of the outward normal $\nu$ with a speed at time zero given by
$$\vph:\partial \ca_0 \to\R,$$ 
we would like to compute $\frac{dh}{dt}$,
or equivalently $\frac{d}{dt}\mathrm{Cap}(\ca_t)$,
at time zero in terms of $\vph$ and $u_0$, where $u_t:\overline{\ca_t}\to [0,1]$ is the  harmonic function that realises the variational problem in \eqref{var_prob} for $\ca_t$,
which is unique if we specify on which evolving boundary component of $\partial\ca_t$ the function $u_t$ is to take the value $0$ rather than $1$. The function $u_t$ varies smoothly in time in the sense that $(x,t)\mapsto u_t(x)$ is a smooth function 
on
$$\Omega:=\{(x,t)\in M\times (-\ep,\ep)\ |\ x\in\overline{\ca_t}\}.$$

We simplify notation by writing $u:=u_0$ and $\ca:=\ca_0$.

We first compute the derivative $\pl{}{t}u_t(x)$ for each $x\in\overline{\ca} $, where in general we must interpret this as a one-sided derivative when $x\in \partial\ca $.
This derivative must give a harmonic function on $\ca $ because each $u_t$ is harmonic, 
so it suffices to compute the derivative at an arbitrary point $x\in \partial\ca $.

Let $\si:(-\ep,\ep)\to M$ be a smooth path so that $\si(t)\in \partial\ca_t$ and so 
that $\dot\si(t)$ is normal to $\partial\ca_t$. Thus we can write 
$\dot\si(0)=\vph \nu$ at $\si(0)$.
Because $u_t(\si(t))$ is constant in $t$ 
(taking the value $0$ or $1$)
we can compute
$$0=\frac{d}{dt}\bigg|_{t=0}u_t(\si(t))=
\vph\, \nu.\grad u
+\pl{u_t}{t}\bigg|_{t=0}
$$
at $\si(0)$. 
Because
\beq
\label{handy_new}
\grad u = \pl{u}{\nu} \nu
\eeq
on the boundary $\partial \ca $,
we deduce that $\dot u:=\pl{u_t}{t}\big|_{t=0}$ is the harmonic function 
on $\ca $ with boundary values given by
$$\dot u = -\vph \pl{u}{\nu}.$$

As a result we can compute the evolution of the energy 
$$E(u_t)=\half\int_{\ca_t} |\grad u_t|^2d\mu,$$ 
with the harmonic function $u_t$ moving in the direction of $\dot{u}$ given as above, as
\beqa
\label{energy_evol_form_old}
\frac{d}{dt}\bigg|_{t=0}E(u_t)
&=
\int_{\ca }\grad \dot{u} . \grad u\, d\mu
+\half\int_{\partial \ca } |\grad u|^2\vph ds\\
&=
\int_{\partial \ca } \dot{u}\,(\nu.\grad u)\, ds
+\half\int_{\partial \ca } |\grad u|^2\vph ds\\
&=
-\int_{\partial \ca } \left|\pl{u}{\nu}\right|^2 \vph\, ds
+\half\int_{\partial \ca } |\grad u|^2\vph ds\\
&=-\half\int_{\partial \ca } |\grad u|^2\vph ds,
\eeqa
where again we used \eqref{handy_new}.
In conclusion, we have computed the evolution of the capacity of an evolving annulus:
\begin{lem}
Suppose $\ga^0_t$ and $\ga^1_t$ are smoothly varying 
disjoint embedded loops within a 
Riemannian surface, for $t\in (-\ep,\ep)$, that bound 
an annulus $\ca_t$ that evolves continuously with respect to the Hausdorff distance.
Suppose that 
$\ga^0_t$ and $\ga^1_t$ move in the outward normal direction with speed $\vph:\partial\ca_0\to\R$ at time $t=0$.
Then the energy $E(u_t)$ of the harmonic function $u_t:\overline{\ca_t}\to[0,1]$ satisfying $u_t\equiv 0$ on $\ga^0_t$ and $u_t\equiv 1$ on $\ga^1_t$ evolves according to
\beq
\label{energy_evol_form}
\frac{d}{dt}\bigg|_{t=0}E(u_t)
=-\half\int_{\partial \ca_0} |\grad u_0|^2\vph\, ds.
\eeq
\end{lem}

\section{Evolving an annulus under curve shortening flow}

We  continue the discussion of the previous section in the special case that
the normal speed of the annulus is given by $\vph=-\ka$, that is, we evolve the boundary components under curve shortening flow.
Our goal is to prove Theorem \ref{main_thm}.
By translating time, we are free to work over time intervals $(-\ep,\ep)$ as in the previous section 
and prove the formula \eqref{logh_ineq} at time $t=0$.

Since the boundary $\partial\ca$ is the level set of a harmonic function $u$, 
the geodesic curvature $\ka$ on $\partial \ca $ can be written $\ka = -\pl{}{\nu}\log |\grad u|$. The most convenient form for this estimate is:

\textbf{Claim:}
\beq
\label{kappa_formula_new}
\half\pl{}{\nu}|\grad u|^2= -\ka |\grad u|^2.
\eeq

\begin{proof}
Let $\tau$ be a vector field that is of unit length and tangent to $\partial\ca$ on $\partial\ca$.
By harmonicity of $u$, we can compute on the level sets of $u$ making up $\partial\ca$ that 
$$\Hess(u)(\nu,\nu)=-\Hess(u)(\tau,\tau)=-\grad_\tau(\grad_\tau u)+\grad_{\grad_\tau \tau}u
=0+\grad_{-\ka \nu} u = -\ka\pl{u}{\nu}.$$
Consequently, using \eqref{handy_new} again, we have 
\beqa
\half\pl{}{\nu}|\grad u|^2 &= \Hess(u)(\nu,\grad u)
= \pl{u}{\nu} \Hess(u)(\nu,\nu)\\
&= -\ka \left|\pl{u}{\nu}\right|^2
=-\ka |\grad u|^2,
\eeqa
which completes the proof of the claim \eqref{kappa_formula_new}.
\end{proof}

We can  use the claim \eqref{kappa_formula_new} in our formula \eqref{energy_evol_form} for the evolution of the energy, still restricted to the case $\vph=-\ka$, to give
\beqa
\frac{d}{dt}\bigg|_{t=0}E(u_t) &=\half\int_{\partial \ca } |\grad u|^2\ka\, ds\\
&=-\frac14 \int_{\partial \ca } \pl{}{\nu} |\grad u|^2 ds\\
&=-\frac14 \int_{\ca } \lap |\grad u|^2 d\mu.
\eeqa
Applying the Bochner formula, using that $u$ is harmonic, then gives
$$\frac{d}{dt}\bigg|_{t=0}E(u_t)=-\half\int_{\ca } (|\Hess(u)|^2+K|\grad u|^2)d\mu,$$
where $K$ is the Gauss curvature of the underlying surface.
If $K\geq K_0$ then this implies that 
$$\frac{d}{dt}\bigg|_{t=0}E(u_t)\leq -K_0 E(u),$$
that is, 
$$\frac{d}{dt}\bigg|_{t=0}\log E(u_t)\leq -K_0.$$
But $\frac{1}{2h}=E(u_t)$, so 
$$\frac{d}{dt}\bigg|_{t=0}\log h\geq K_0,$$
as required.
In the case of equality we must have $\Hess(u)\equiv 0$. This forces $\grad u$ to be a Killing field of constant length.
In our special case this implies that every 
level set of $u$ is isometric to a circle of fixed length
and the domain $\ca $ of $u$ must be isometric to a flat cylinder of finite length. This completes the proof of Theorem \ref{main_thm}.

\brmk
The proof given above simplifies an earlier presentation in which we worked directly with the biholomorphic maps from the annuli to cylinders $(0,h)\times S^1$. That approach 
can be found in \cite{arjun_thesis}.
\ermk

\vskip .2cm

\noindent
AS: 
\url{https://sites.google.com/view/sobnack}

\noindent
{\sc Department of Mathematics, The University of Texas at Austin, TX 78712, USA.} 

\noindent
PT: 
\url{https://warwick.ac.uk/fac/sci/maths/people/staff/peter_topping}

\noindent
{\sc Mathematics Institute, University of Warwick, Coventry,
CV4 7AL, UK.}

\end{document}